\documentclass[12pt,thmsa]{article}
\usepackage{amsfonts, amsmath}

\newtheorem{theorem}{Theorem}

\newtheorem{proposition}{Proposition}
\newtheorem{lemma}{Lemma}

\newcommand{\p}{\Bbb{P}}

\newcommand{\e}{\Bbb{E}}

\newcommand{\ind}{\mbox{\rm 1\hspace{-0.04in}I}}

\newcommand{\R}{\mbox{\rm I\hspace{-0.02in}R}}

\newcommand{\ud}{\mathrm{d}}

\def\QED{\hfill\vrule height 1.5ex width 1.4ex depth -.1ex \vskip20pt}

\begin{document}
\hspace*{-0.5in} {\footnotesize This version March 19, 2007.}
\vspace*{0.9in}
\begin{center}
{\LARGE  On the rate of growth of L\'evy processes with no positive
jumps conditioned to
stay positive .\vspace*{0.4in}}\\
{\large J.C. Pardo\footnote{Research supported by a grant from
CONACYT (Mexico).}\vspace*{0.2in}}
\end{center}

\noindent $^{1}$  {\footnotesize Laboratoire de Probabilit\'es et
Mod\`eles Al\'eatoires, Universit\'e Pierre et Marie Curie, 4,
Place Jussieu - 75252 {\sc Paris Cedex 05.} E-mail:
pardomil@ccr.jussieu.fr}\\
\author{J.C. Pardo.}\\

\noindent {\it Abstract}: {\footnotesize In this article, we study
the asymptotic behaviour of L\'evy processes with no positive jumps
conditioned to stay positive. We establish integral tests for the
lower envelope at $0$ and at $+\infty$ and an analogue of
Khintchin's law of the iterated logarithm at $0$ and $+\infty$, for the upper envelope.} \\


\noindent
{\it A.M.S. Classification}: {\footnotesize  60 G 17, 60 G 51.}\\
\section{Introduction and main results.}
L\'evy processes conditioned to stay positive have been introduced
by Bertoin at the beginning of the nineties (see \cite{be2},
\cite{be3}, \cite{be4}). These first studies are devoted to the
special case where the L\'evy processes are spectrally one-sided. In
a recent work Chaumont
and Doney \cite{chd} studied more general cases.\\
In this article, we are interested in the case when L\'evy
processes have no positive jumps (or spectrally negative L\'evy processes). This case has been deeply studied
by Bertoin in Chapter VII of \cite{be}. This
will be our basic reference.\\
The aim of this note is to describe the lower and upper envelope at
$0$ and at $+\infty$ of  L\'evy processes with no positive jumps
conditioned to stay positive throughout integral tests and laws
of the iterated logarithm.\\
For our purpose, we will first introduce some important properties of L\'evy
processes with no positive jumps and then we will define the
``conditioned" version.\\
Let $\mathcal{D}$ denote the Skorokhod's space of c\`adl\`ag paths
with real values and defined on the positive real half-line
$[0,\infty)$ and $\p$ a probability measure defined on $\mathcal{D}$
under which $\xi$ will  be a real-valued L\'evy process with no
positive jumps, that is its L\'evy measure has support in the
negative real-half line.\\
From the general theory of L\'evy processes (see Bertoin \cite{be}
or Sato \cite{Sa} for background), we know that $\xi$ has finite
exponential moments of arbitrary positive order. In particular $\xi$
satisfies
\[
\e\Big(\exp\big\{\lambda\xi_t\big\}\Big)=\exp\big\{t\psi(\lambda)
\big\},\qquad \lambda, t\geq 0,
\]
where
\[
\psi(\lambda)=a\lambda+\frac12 \sigma^2 \lambda^2+\int_{(-\infty,
0)}\big(e^{\lambda x}-1-\lambda x\ind_{\{x>-1\}}\big)\Pi(\ud x),
\quad \lambda \geq 0,
\]
$a\in \R, \sigma\geq 0$ and $\Pi$ is a measure that satisfies
\[
\int_{(-\infty, 0)}(1\land x^2)\Pi(\ud x).
\]
The measure $\Pi$ is well-known as the L\'evy measure of the process $\xi$.\\
According to Bertoin \cite{be}, the mapping $\psi:[0, \infty)\to(-\infty, \infty)$ is convex and
ultimately increasing. We denote its right-inverse on $[0,
\infty)$ by $\Phi$. Let us introduce the first passage time of $\xi$ by
\[
T_x=\inf\big\{s:\xi_s\ge x\big\}\qquad
\textrm{for}\qquad x\geq0.
\]
From Theorem VII.1 in \cite{be}, we know that
the  process $T=(T_x, x\geq 0)$ is a subordinator, killed at an independent
exponential time if $\xi$ drifts towards $-\infty$. The Laplace exponent of $T$ is given by $\Phi$,
\[
\e\Big(\exp\big\{-\lambda T_x\big\}\Big)=\exp\big\{-x\Phi(\lambda)
\big\},\qquad \lambda, t\geq 0.
\]
According to
Bertoin \cite{be} Chapter III, we know that
\[
\Phi(\lambda)=\mathrm{k} +\mathrm{d} \lambda
+\int_{(0,\infty)}(1-e^{-\lambda x}) \nu(\ud x),\qquad \lambda\geq
0,
\]
where $\mathrm{k}$ is the killing rate, $\mathrm{d}$ is the drift
coefficient and $\nu$ is the L\'evy measure of the subordinator $T$
which fulfils the following condition,
\[
\int_{(0, \infty)}(1\land x)\nu(\ud x)<\infty.
\]
It is sometimes convenient to perform an integration by parts and
rewrite $\Phi$ as
\[
\frac{\Phi(\lambda)}{\lambda}=\mathrm{d}+\int_{0}^{\infty}
e^{-\lambda x}\bar{\nu}(x)\ud x, \qquad \textrm{with }\quad
\bar{\nu}(x)=\mathrm{k} + \nu \big((x,\infty)\big).
\]
Note that the killing rate and the drift coefficient are given by
\[
\mathrm{k}=\Phi(0)\qquad\textrm{ and }\qquad
\mathrm{d}=\lim_{\lambda\to +\infty}\frac{\Phi(\lambda)}{\lambda}.
\]
In particular, the life time $\zeta$ has an exponential distribution
with parameter $k\geq 0$ ($\zeta=+\infty$ for $\mathrm{k}=0$).
Hence, it is not difficult to deduce that  $\xi$ drifts towards
$-\infty$ if and only if $\Phi(0)>0$.\\
In order to study the case when $\xi$ drifts towards $-\infty$, we
will define the following probability measure,
\[
\p^{\natural}(A)=\e\Big(\exp\big\{\Phi(0)\xi_{t}\big\} \ind_A\Big),
\qquad A \in \mathcal{F}_{t},
\]
where $\mathcal{F}_{t}$ is the $\p$-complete sigma-field generated
by $(\xi_{s}, s\leq t)$. Note that  under $\p^{\natural}$, the
process $\xi$ is still a L\'evy process with no positive jumps which
drifts towards $+\infty$ and its Laplace exponent is defined by
$\psi^{\natural}(\lambda)=\psi(\Phi(0)+\lambda)$, $\lambda\geq 0$.
Moreover the first passage process $T$ is still a subordinator with
Laplace exponent $\Phi^{\natural}(\lambda)=\Phi(\lambda)-\Phi(0)$.\\
Since $\xi$ has no positive jumps, we can solve the so-called
two-sided-exit problem in  explicit form. This  problem consists
in determining the probability that $\xi$ makes its first exit from
an interval $[-x, y]$ ($x,y>0$) at the upper boundary point. More
precisely,
\[
\p\bigg(\inf_{0\leq t\leq T_y}\xi_t\geq
-x\bigg)=\frac{W(x)}{W(x+y)},
\]
where $W:[0, \infty)\to[0, \infty)$ is the unique absolutely
continuous increasing function with Laplace transform
\begin{equation}\label{scfunc}
\int_{0}^{\infty}e^{-\lambda x}W(x)\ud x=\frac{1}{\psi(\lambda)},
\quad \textrm{for} \quad \lambda>\Phi(0).
\end{equation}
The function $W$ is well-known as the scale function and is necessary for
the definition of  L\'evy processes conditioned to stay
positive.\\
Using the Doob's theory of $h$-transforms, we construct a new Markov
process by an $h$-transform  of the law of the L\'evy process killed
at time $R=\inf\{t\geq 0: \xi_t<0\}$ with the harmonic function $W$
(see for instance Chapter VII in Bertoin \cite{be}, Chaumont
\cite{ch} or Chaumont and Doney \cite{chd}), and its semigroup is
given by
\[
\p^{\uparrow}_x(\xi_t\in\ud y)=\frac{W(y)}{W(x)}\p_x(\xi_t\in \ud
y, t<R)\quad \textrm{ for }\quad x>0,
\]
where $\p_x$ denotes the law of $\xi$ starting from $x> 0$. Under
$\p^{\uparrow}_x$, $\xi$ is a process taking values in $[0,
\infty)$. It will be referred to as the L\'evy process started at
$x$ and conditioned to stay positive.\\
It is important to note that when $\xi$ drifts towards $-\infty$,
 we have that $\p^{\uparrow}_{x}=\p^{\natural \uparrow}_{x}$, for all $x\geq
 0$. Hence the study of this case is reduced to the study of the
 processes which drift towards $+\infty$.\\
Bertoin proved in \cite{be} the existence of a measure
$\p^{\uparrow}_0$ under which the process starts at $0$ and stays
positive. In fact, the author in \cite{be} proved that the probability measures
$\p^{\uparrow}_x$ converge as $x$ goes to $0+$ in the sense of
finite-dimensional distributions to
$\p^{\uparrow}_0:=\p^{\uparrow}$ and noted that this convergence
also holds in the sense of Skorokhod. Several authors have studied
this convergence, the most recent work is due to Chaumont and
Doney \cite{chd}. In \cite{chd}, the authors proved that this
convergence holds in the sense of Skorokhod
under general hypothesis.\\
Chaumont and Doney \cite{chd} also give a path decomposition of the
process $(\xi, \p^{\uparrow}_x)$ at the time of its minimum. In
particular, if $m$ is the time at which $\xi$, under
$\p^{\uparrow}_x$, attains its minimum, we have that under
$\p^{\uparrow}_x$ the pre-infimum process, $(\xi_t, 0\leq t<m)$, and
the post-infimum process, $(\xi_{t+m}-\xi_{m}, t\geq 0)$, are
independent and the later has by  law  $\p^{\uparrow}$. Moreover,
the process $(\xi, \p^{\uparrow}_x)$ reaches its absolute minimum
once only and its law is given by
\begin{equation}\label{distabsmin}
 \p^{\uparrow}_x( \xi_{m}\geq y)=\frac{W(x-y)}{W(x)}\ind_{\{y\leq
 x\}}.
\end{equation}
Recently, Chaumont and Pardo \cite{CP} studied the lower envelope of
positive self-similar Markov processes (or pssMp for short). In particular, the
authors obtained integral tests at $0$ and at $+\infty$, for the
lower envelope of stable L\'evy processes with no positive jumps
conditioned to stay positive  and with index $\alpha \in (1, 2]$;
such processes are well-known examples of pssMp. More precisely,
when $X^{(x)}$ is the stable L\'evy process with no positive jumps
conditioned to stay positive and with starting point $x\ge 0$, we
have the following integral test for the lower envelope at $0$ and
at $+\infty$: let $f$ be an increasing function, then
\[
\p(X^{(0)}_t<f(t),\mbox{i.o., as
$t\rightarrow0$})=\left\{\begin{array}{l}0\\1\end{array}\right.\,,
\textrm{ according as }\int_{0+}\left(\frac
{f(t)}t\right)^{1/\alpha}\,\frac{dt}t\;
\left\{\begin{array}{l}<\infty\\=\infty\end{array}\right.\,.\] Let
$f$ be an increasing function, then for all $x\ge0$,
\[\p(X^{(x)}_t<f(t),\mbox{i.o., as
$t\rightarrow\infty$})=\left\{\begin{array}{l}0\\1\end{array}\right.\,,
\textrm{ according as } \int^{+\infty}\left(\frac
{f(t)}t\right)^{1/\alpha}\,\frac{dt}t\;
\left\{\begin{array}{l}<\infty\\=\infty\end{array}\right.\,.\]
Later, Pardo \cite{Pa1} studied the upper envelope of pssMp  at $0$
and at $+\infty$. In particular, the author noted that $X^{(x)}$, the stable L\'evy 
process with no positive jumps conditioned to stay positive,
satisfies the following law of the iterated logarithm
\[
\limsup_{t\to 0}\frac{X^{(0)}_t}{t^{1/\alpha}(\log|\log
t|)^{1-1/\alpha}}=c(\alpha), \qquad \p^{\uparrow}-a.s.,
\]
where $c(\alpha)$ is a positive constant which depends on $\alpha$.\\
Moreover, Pardo \cite{Pa1} also established that
\begin{equation*}
\limsup_{t\to 0}\frac{J^{(0)}_{t}}{t^{1/\alpha}(\log|\log
t|)^{1-1/\alpha}}=c(\alpha), \qquad \p^{\uparrow}-a.s.,
\end{equation*}
\begin{equation*}
\limsup_{t\to 0}\frac{X^{(0)}_t-J^{(0)}_{t}}{t^{1/\alpha}(\log|\log
t|)^{1-1/\alpha}}=c(\alpha), \qquad \p^{\uparrow}-a.s.,
\end{equation*}
where $J^{(x)}$ is the future infimum process of $X^{(x)}$, i.e.
$J^{(x)}_t=\inf_{s\geq t}X^{(x)}_s$. It is important to note that
the above laws of the iterated
logarithm have been also obtained at $+\infty$, for any starting point $x\geq 0$.\\
Bertoin \cite{be1} studied the local rate of growth of L\'evy
processes with no positive jumps. In \cite{be1}, the author  noted that
the sample path behaviour of a L\'evy process with no positive jumps $\xi$, immediately
after a local minimum is the same as that of its conditioned version
$(\xi, \p^{\uparrow})$ at the origin. The main  result in \cite{be1} gives us 
a remarkable law of the iterated logarithm at an
instant when the path attains a local minimum on the interval $[0,
1]$. We will recall this result for L\'evy processes conditioned to
stay positive.\\
With a misuse of notation, we will denote by $\Phi$ and $\psi$ for
the functions $\Phi^{\natural}$ and $\psi^{\natural}$, respectively;
when we are in the case where the process $\xi$ drifts towards
$-\infty$.
\begin{theorem}[Bertoin, \cite{be1}]
There is a finite positive constant $k$ such that
\[
\limsup_{t\to 0}\frac{\xi_t\Phi(t^{-1}\log|\log t|)}{\log|\log
t|}=k, \qquad \p^{\uparrow}-a.s.
\]
\end{theorem}
It is important to note that  the constant found by Bertoin
satisfies $k\in[c,c+\gamma]$, where $c$ is the same constant found
below in Theorem 3 and  $\gamma\ge 6$.

Bertoin presented in \cite{be1} a ``geometric" proof using some path
transformations that connect $\xi$ and its conditioned version
$(\xi, \p^{\uparrow})$. Here, we will present standard arguments
involving probability
tail estimates.\\
Our main results requires the following  hypothesis, for
all $\beta>1$
\begin{equation*}
\textrm{(H1)}\quad\liminf_{x\to 0}\frac{\psi(x)}{\psi(\beta x)}>0
\quad\textrm{and}\quad \textrm{(H2)}\quad\liminf_{x\to
+\infty}\frac{\psi(x)}{\psi(\beta x)}>0.
\end{equation*}
\begin{theorem}
Let us suppose that (H2) is satisfaied, then there is a
positive constant $k$ such that,
\[
\limsup_{t\to 0}\frac{\xi_t\Phi(t^{-1}\log|\log t|)}{\log|\log
t|}=k, \qquad \p^{\uparrow}-a.s.
\]
Moreover, if condition (H1) is satisfied, then 
\[
\limsup_{t\to +\infty}\frac{\xi_t\Phi(t^{-1}\log\log t)}{\log|\log
t|}=k, \qquad \p^{\uparrow}-a.s.
\]
\end{theorem}
Note that with our arguments, we found that $k\in [c,c\eta]$
where of course $\eta\ge 1$ and $c\eta> 3$. We also remark that in
particular (H1) and (H2) are satisfied under the assumption that
$\psi$ is regularly varying at $0$ and at $\infty$ with index
$\alpha>1$. Under this assumption
\[
k=c=(1/\alpha)^{-1/\alpha}(1-1/\alpha)^{\frac{1-\alpha}{\alpha}}.
\]
The next result gives us a law of the iterated logarithm at $0$
and at $+\infty$ of the future infimum of $(\xi, \p^{\uparrow})$.
This result extends the result of Pardo \cite{Pa1} for the stable case.
\begin{theorem}Let $J$ denote the future infimum process of $\xi$,
defined by $J_t=\inf_{s\geq t}\xi_s$, then there is a positive
constant $c$ such that,
\[
\limsup_{t\to 0}\frac{J_t\Phi(t^{-1}\log|\log t|)}{\log|\log
t|}=c, \qquad \p^{\uparrow}-a.s.,
\]
and
\[
\limsup_{t\to +\infty}\frac{J_t\Phi(t^{-1}\log\log t)}{\log|\log
t|}=c, \qquad \p^{\uparrow}-a.s.
\]
If we assume that (H2) is
satisfied, then 
\[
\limsup_{t\to 0}\frac{J_t\Phi(t^{-1}\log|\log t|)}{\log|\log t|}=k,
\qquad \p^{\uparrow}-a.s.,
\]
i.e. that $c=k$.\\
Moreover if (H1) is satisfied, then
\[
\limsup_{t\to +\infty}\frac{J_t\Phi(t^{-1}\log\log t)}{\log|\log
t|}=k, \qquad \p^{\uparrow}-a.s.,
\]
i.e. that $c=k$.
\end{theorem}
We now turn our attention to the L\'evy process conditioned to
stay positive reflected at its future infimum. The following
theorem extends the law of the iterated logarithm of Pardo \cite{Pa1} for the stable case.\\
Let us suppose that for all $\beta <1$
\begin{equation*}
\textrm{(H3)} \quad\limsup_{x\to 0}\frac{W(\beta x)}{W(x)}<1 \quad \textrm{and}\quad
\textrm{(H4)}\quad\limsup_{x\to +\infty}\frac{W(\beta x)}{W(x)}<1.
\end{equation*}
\begin{theorem}
 Under the hypothesis (H2) and (H3), we have that
\[
\limsup_{t\to 0}\frac{\big(\xi_t - J_t\big)\Phi(t^{-1}\log|\log
t|)}{\log|\log t|}=k, \qquad \p^{\uparrow}-a.s.
\]
Moreover if (H1) and (H4) are satisfied, then
\[
\limsup_{t\to +\infty}\frac{\big(\xi_t - J_t\big)\Phi(t^{-1}\log\log
t)}{\log|\log t|}=k, \qquad \p^{\uparrow}-a.s.
\]
\end{theorem}
Again conditions (H3) and (H4) are  satisfied when  $\psi$ is regularly varying
at $0$ and at $\infty$ with index $\alpha>1$.\\
The lower envelope of $(\xi, \p^{\uparrow})$ at $0$ and at $\infty$ is determined as follows,
\begin{theorem}  If $\xi$ has unbounded variation and $f:[0, \infty)\to [0, \infty)$ is an increasing function such that
$t\to f(t)/t$ decreases, one has
\[
\liminf_{t\to 0}\frac{\xi_t}{f(t)}=0\quad \p^{\uparrow}\textrm{-a.s.} \quad\textrm{if and only if}\quad \int_{0}f(x)\nu(\ud x)=\infty.
\]
Moreover,
\[
\textrm{if }\qquad\int_{0}f(x)\nu(\ud x)<\infty\qquad\textrm{ then }\qquad
\lim_{t\to 0}\frac{\xi_t}{f(t)}=\infty \qquad \p^{\uparrow}\textrm{-a.s.}
\]
The lower envelope at $+\infty$ is determined as follows: if $\xi$ oscillates or drifts towards $-\infty$ and the function $f:[0, \infty)\to [0, \infty)$ is
 increasing  such that $t\to f(t)/t$ decreases, one has
\[
\liminf_{t\to +\infty }\frac{\xi_t}{f(t)}=0\quad \p^{\uparrow}\textrm{-a.s.} \quad\textrm{if and only if}\quad \int^{+\infty}f(x)\nu(\ud x)=\infty.
\]
Moreover,
\[
\textrm{if }\qquad\int^{+\infty}f(x)\nu(\ud x)<\infty\qquad \textrm{ then }\qquad
\lim_{t\to +\infty}\frac{\xi_t}{f(t)}=\infty \qquad \p^{\uparrow}\textrm{-a.s.}
\]
\end{theorem}
The following result describes the lower envelope of the future infimum of $(\xi, \p^{\uparrow})$. In fact, we will
deduce that $(\xi, \p^{\uparrow})$ and its future infimum have the same lower functions.
\begin{theorem}
$\mathrm{(i)}$ If $\xi$ has bounded variation, one has
\[
\lim_{t\to 0}\frac{J_t}{t}=\frac{1}{\mathrm{d}}\qquad \p^{\uparrow}\textrm{-a.s.}
\]
$\mathrm{(ii)}$ If $\xi$ has unbounded variation and $f:[0, \infty)\to [0, \infty)$ is an increasing function such that
$t\to f(t)/t$ decreases, one has
\[
\liminf_{t\to 0}\frac{J_t}{f(t)}=0\quad \p^{\uparrow}\textrm{-a.s.} \quad\textrm{if and only if}\quad \int_{0}f(x)\nu(\ud x)=\infty.
\]
Moreover, 
\[
\textrm{if }\qquad \int_{0}f(x)\nu(\ud x)<\infty \qquad \textrm{ then }\qquad
\lim_{t\to 0}\frac{J_t}{f(t)}=\infty \qquad
\p^{\uparrow}\textrm{-a.s.}
\]
$\mathrm{(iii)}$ If $\xi$ drifts to $+\infty$ one has
\[
\lim_{t\to +\infty}\frac{J_t}{t}=\frac{1}{\e(T_1)}\qquad
\p^{\uparrow}\textrm{-a.s.}
\]
$\mathrm{(iv)}$ If $\xi$ oscillates or drifts to $-\infty$ and $f:[0, \infty)\to [0, \infty)$ is an increasing function such that the mapping
$t\to f(t)/t$ decreases, one has
\[
\liminf_{t\to +\infty }\frac{J_t}{f(t)}=0\quad \p^{\uparrow}\textrm{-a.s.} \quad\textrm{if and only if}\quad \int^{+\infty}f(x)\nu(\ud x)=\infty.
\]
Moreover, 
\[
\textrm{if }\qquad\int^{+\infty}f(x)\nu(\ud x)<\infty\qquad \textrm{ then }\qquad
\lim_{t\to +\infty}\frac{J_t}{f(t)}=\infty \qquad
\p^{\uparrow}\textrm{-a.s.}
\]
\end{theorem}
Moreover, we have the following integral test for the lower functions in terms of $\Phi$.
\begin{proposition}$\mathrm{(i)}$ Let $f:[0, \infty)\to [0, \infty)$ be an increasing function.
\[
\textrm{If}\qquad\int_0 x^{-1}f(x)\Phi(1/x)\ud x<\infty\qquad \textrm{ then }\qquad
\lim_{t\to 0}\frac{J_t}{f(t)}=\lim_{t\to 0}\frac{\xi_t}{f(t)}=\infty \qquad \p^{\uparrow}\textrm{-a.s.}
\]
$\mathrm{(ii)}$ Let $f:[0, \infty)\to [0, \infty)$ be an increasing function.
\[
\textrm{If}\qquad\int^{+\infty} x^{-1}f(x)\Phi(1/x)\ud x<\infty\quad \textrm{ then }\quad
\lim_{t\to +\infty}\frac{J_t}{f(t)}=\lim_{t\to +\infty}\frac{\xi_t}{f(t)}=\infty \quad \p^{\uparrow}\textrm{-a.s.}
\]
\end{proposition}
The rest of this note consists of two sections, which are devoted to
the following topics: Section 2 provides asymptotic results for the
first and the last passage times of the process $(\xi,
\p^{\uparrow})$. In Section 3, we will prove the results presented
above.
\section{First and  last passage times.}
Let us recall the definition of the first and last passage time of
$(\xi,\p^{\uparrow})$ or $\xi^{\uparrow}$ to simplify our
notation,
\[
T^{\uparrow}_{x}=\inf\Big\{t\geq 0 : \xi^{\uparrow}_t\geq
x\Big\}\quad\textrm{and}\quad U^{\uparrow}_{x}=\sup\Big\{t\geq 0 :
\xi^{\uparrow}_t\leq x\Big\}\quad\textrm{for}\quad x\geq 0.
\]
From Theorem VII.18 in \cite{be}, we know that $(\xi_t, 0\leq t\leq
T_x)$, the L\'evy process killed at its first passage time above
$x$, under $\p^{\natural}$, has the same law as the L\'evy process
conditioned to stay positive time-reversed at its last passage time
below $x$, $(x-\xi^{\uparrow}_{(U^{\uparrow}_x-t)^-}, 0\leq t\leq
U^{\uparrow}_x)$. In particular,  we deduce that $U^{\uparrow}_x$
has the same law as $T_x$ and that $U^{\uparrow}=(U^{\uparrow}_x,
x\geq 0)$ is a subordinator with Laplace exponent $\Phi(\lambda)$
and therefore we obtain that the process $\xi^{\uparrow}$ drifts
towards $+\infty$.\\
There exist a huge variety of results on  the upper envelope of
subordinators. Fristedt and Pruitt \cite{FP} proved a general law of
the iterated logarithm which is valid for a wide class of
subordinators. The sharper result on the lower envelope for
subordinators is due to Pruitt \cite{Pr}. In his main result, he
gave an important integral test.\\
Bertoin \cite{be} presents a more precise law of the iterated
logarithm of subordinators than the result obtained by Fristedt and
Pruitt but for a more restrictive class of subordinators. In his
result, Bertoin supposes that $\psi$ is regularly varying at
$+\infty$ with index $\alpha>1$ (see Theorems III.11 and III.14 in
\cite{be}). In particular, we have the following lemma.
\begin{lemma}
The last passage time process $U^{\uparrow}$ under the assumption
that $\psi$ is regularly varying at $+\infty$ with index $\alpha>1$,
satisfies
\[
\liminf_{x\to 0}\frac{U^{\uparrow}_x \psi\big(x^{-1}\log|\log
x|\big) }{\log|\log
x|}=\frac{1}{\alpha}\left(1-\frac{1}{\alpha}\right)^{\alpha-1},\qquad
\text{almost surely,}
\]
and for large times, if we suppose that $\psi$ is regularly varying
at $0$ with index $\alpha>1$, then
\[
\liminf_{x\to +\infty}\frac{U^{\uparrow}_x \psi\big(x^{-1}\log\log
x\big)}{\log\log
x}=\frac{1}{\alpha}\left(1-\frac{1}{\alpha}\right)^{\alpha-1},\qquad
\text{almost surely}.
\]
\end{lemma}
Now, we turn our attention to the first passage time process. Note
that due to the absence of positive jumps, for all $x\geq 0$,
$\xi^{\uparrow}_{T^{\uparrow}_x}=x$, a.s. Hence from the strong
Markov property, we have that $T^{\uparrow}=(T^{\uparrow}_x, x \geq
0)$ is an increasing process with independent increments but not
stationary.\\
Here we will use the results presented in Bertoin \cite{be} and
Lemma 1 to obtain the following law of the iterated logarithm for
the first and last passage time of $\xi^{\uparrow}$.
\begin{proposition}Suppose that $\psi$ is regularly varying at
$+\infty$ with index $\alpha>1$. Then the first passage time process
satisfies the following law of the iterated logarithm,
\[
\liminf_{x\to 0}\frac{T^{\uparrow}_x \psi\big(x^{-1}\log|\log
x|\big)}{\log|\log
x|}=\frac{1}{\alpha}\left(1-\frac{1}{\alpha}\right)^{\alpha-1},\qquad
\text{almost surely,}
\]
and for large times, if we suppose that $\psi$ is regularly
varying at $0$ with index $\alpha>1$, we have
\[
\liminf_{x\to +\infty}\frac{T^{\uparrow}_x \psi\big(x^{-1}\log\log
x\big)}{\log\log
x}=\frac{1}{\alpha}\left(1-\frac{1}{\alpha}\right)^{\alpha-1},\qquad
\text{almost surely}.
\]
\end{proposition}
\textit{Proof of Proposition 2:} We will only prove the result for small times since
the proof for large times is very similar. For all $x\geq 0$, we see
that $T^{\uparrow}_x\leq U^{\uparrow}_x$, then from Lemma 1 we
obtain the upper bound
\[
\liminf_{x\to 0}\frac{T^{\uparrow}_x \psi\big(x^{-1}\log|\log
x|\big)}{\log|\log x|}\leq\liminf_{x\to
0}\frac{U^{\uparrow}_x\psi\big(x^{-1}\log|\log x|\big)}{\log|\log
x|}=\frac{1}{\alpha}\left(1-\frac{1}{\alpha}\right)^{\alpha-1}.
\]
Next, we prove the lower bound. With this purpose we establish the
following lemma.
\begin{lemma}
Assume that $\psi$ is regularly varying at $+\infty$ with index
$\alpha>1$. Then for every constant $c_1>0$, we have
\[
-\log \p\Big(T^{\uparrow}_x\leq cg(x)\Big)\sim
\left(1-\frac{1}{\alpha}\right)\left(\frac{1}{c_1\alpha}\right)^{1/(\alpha-1)}\log|\log
x|, \qquad\textrm{as}\quad x\to 0,
\]
where
\[
g(t)=\frac{\log|\log x|}{\psi\big(x^{-1}\log|\log x|\big)}.
\]
\end{lemma}
\textit{Proof of Lemma 2:} We know that $U^{\uparrow}$ is a subordinator and that $\psi$
is the inverse of the function $\Phi$, then from Lemma
III.12 in Bertoin \cite{be}, we see that
\[
-\log \p\Big(U^{\uparrow}_x\leq c_1 g(x)\Big)\sim
\left(1-\frac{1}{\alpha}\right)\left(\frac{1}{c_1\alpha}\right)^{1/(\alpha-1)}\log|\log
x|, \qquad\textrm{as}\quad x\to 0.
\]
Then the upper bound is clear since for all $x>0$ we have that
$T^{\uparrow}_x\leq U^{\uparrow}_x$.\\
For the lower bound, let us first define the supremum process
$S=(S_t, t \geq 0)$ by $S_t=\sup_{0\leq s\leq t}\xi_s$. Next, we
fix $\epsilon>0$, then by the Markov property
\begin{equation}\label{deslogr}
\begin{split}
\p^{\uparrow}\big(J_{c_1
g(x)}>(1-\epsilon)x\big)&\geq\p^{\uparrow}\big(S_{c_1 g(x)}>x,
J_{c_1 g(x)}>(1-\epsilon)x\big)\\
&=\int_{0}^{c_1 g(x)}\p\Big(T^{\uparrow}_x\in \ud
t\Big)\p^{\uparrow}_x\big(J_{c_1g(x)-t}>(1-\epsilon)x\big)\\
&\geq\p\Big(T^{\uparrow}_x< c_1
g(x)\Big)\p^{\uparrow}_x\big(J_{0}>(1-\epsilon)x\big).
\end{split}
\end{equation}
From the definition of the future infimum process, it is clear
that $J_{0}$ is the absolute minimum of $(\xi, \p^{\uparrow}_x)$
then by (\ref{distabsmin})
\[
\p^{\uparrow}_x\big(J_{0}>(1-\epsilon)x\big)=\frac{W(\epsilon
x)}{W(x)}.
\]
On the other hand, from $(\ref{scfunc})$ and applying the
Tauberian and  Monotone density theorems (see for instance Bertoin
\cite{be} or Bingham et al \cite{Bing}) we deduce that
\[
W(x)\sim\frac{\alpha}{\Gamma(1+\alpha)}\frac{1}{x\psi(1/x)}
\qquad\textrm{as}\quad{x\to 0,}
\]
hence,
\begin{equation}\label{min}
\p^{\uparrow}_x\big(J_{0}>(1-\epsilon)x\big)\to
\epsilon^{(\alpha-1)}\quad\textrm{as }\quad x\to 0.
\end{equation}
Now, since the last passage time process is the right inverse of the
future infimum process, we have that
\[
\p^{\uparrow}\big(J_{c_1
g(x)}>(1-\epsilon)x\big)=\p\Big(U^{\uparrow}_{(1-\epsilon)x}<c_1
g(x)\Big),
\]
and applying Chebyshev's inequality, we have that for every
$\lambda>0$
\[
\p\Big(U^{\uparrow}_{(1-\epsilon)x}<c_1 g(x)\Big)\leq
\exp\Big\{\lambda c_1 g(x)-(1-\epsilon)x\Phi(\lambda)\Big\},
\]
and thus
\[
-\log \p\Big(U^{\uparrow}_{(1-\epsilon)x}<c_1 g(x)\Big)\geq -\lambda
c_1g(x)+(1-\epsilon)x\Phi(\lambda).
\]
Next, we choose $\lambda=\lambda(x)$ such that
$(1-\epsilon)x\Phi(\lambda)=K \log|\log x|$ for some positive
constant $K$, that will be specified later on, then $\lambda=\psi(K
(1-\epsilon)^{-1}x^{-1}\log|\log x|)$. Since $\psi$ is regularly
varying at $\infty$ with index $\alpha$, we see that
\[
\lambda=\lambda(x)\sim
K^{\alpha}(1-\epsilon)^{-\alpha}\psi(x^{-1}\log|\log x|).
\]
This implies
\[
-\lambda c_1 g(x)+(1-\epsilon)x\Phi(\lambda)\sim (K-c_1 K^{\alpha}
(1-\epsilon)^{-\alpha})\log |\log x| \qquad (x \to 0).
\]
We now choose $K$ in such way that $K-c_1 K^{\alpha}
(1-\epsilon)^{-\alpha}$ is maximal, that is
\[
K=(1-\epsilon)^{\alpha/(\alpha-1)}\left(\frac{1}{c_1\alpha}\right)^{1/(\alpha-1)},
\]\
and
\[
K-c_1K^{\alpha}
(1-\epsilon)^{-\alpha}=(1-\epsilon)^{\alpha/(\alpha-1)}\left(\frac{1}{c_1\alpha}\right)^{1/(\alpha-1)}\left(1-\frac{1}{\alpha}\right).
\]
In conclusion, we have established that
\begin{equation}\label{lblast}
(1-\epsilon)^{\alpha/(\alpha
-1)}\left(1-\frac{1}{\alpha}\right)\left(\frac{1}{c_1\alpha}\right)^{1/(\alpha-1)}\leq\liminf_{x\to
0}\frac{-\log\p\Big(U^{\uparrow}_{(1-\epsilon)x}\leq c_1
g(x)\Big)}{\log|\log x|}
\end{equation}
Hence from the inequality (\ref{deslogr}) and (\ref{min}) and
(\ref{lblast}), we deduce
\[
(1-\epsilon)^{\alpha/(\alpha-1)}\left(1-\frac{1}{\alpha}\right)\left(\frac{1}{c_1\alpha}\right)^{1/(\alpha-1)}\leq\liminf_{x\to
0}\frac{-\log\p\Big(T^{\uparrow}_x\leq c_1 g(x)\Big)}{\log|\log x|},
\]
and since $\epsilon$ can be chosen arbitrarily small, the lemma is
proved.\QED \noindent Now we can prove the lower bound of the law of
the iterated logarithm for $T^{\uparrow}$. Let $(x_n)$ be a
decreasing sequence of positive real numbers which converges to $0$
and let us define the event
$A_n=\{T^{\uparrow}_{x_{n+1}}<c_1g(x_n)\}$. Now, we choose
$x_n=r^n$, for $r<1$. From the first Borel-Cantelli's Lemma, if
$\sum_{n}\p(A_n)<\infty$, it follows
\[
T^{\uparrow}_{r^{n+1}}\geq c_1g(r^n)\qquad\textrm{almost surely,}
\]
for all large $n$. Since the function $g$ and the process
$T^{\uparrow}$ are increasing, we have
\[
T^{\uparrow}_x\geq cg(x)\qquad\textrm{for }\quad r^{n+1}\leq x\leq
r^n .
\]
Then, it is enough to prove that $\sum_{n}\p(A_n)<\infty$. In this
direction, we take
\[
0<c_1<c'<\left(\frac{1}{\alpha}\right)^{\alpha}(\alpha-1)^{\alpha
-1}.
\]
Since $\psi$ is
regularly varying and we can chose $r$ close enough to $1$, we see
that for $n_0$ sufficiently large
\[
\sum_{n\geq n_0}\p(A_n)\leq \sum_{n\geq
n_0}\p\Big(T^{\uparrow}_{r^{n+1}}<c'g(r^{n+2})\Big) \leq
\int_{0}^{r^{n_0+2}}\p\Big(T^{\uparrow}_x\leq
c'g\big(x\big)\Big)\frac{\ud x}{x},
\]
and from Lemma 2 this last integral is finite since
\[
\left(1-\frac{1}{\alpha}\right)\left(\frac{1}{c'\alpha}\right)^{1/(\alpha-1)}>1,
\]
with this we finish the proof.\QED
There also exist a huge variety of results for
the upper envelope of subordinators, see for instance Chapter III of
Bertoin \cite{be}. Here, we will state with out proofs the main
results for the upper envelope of $U^{\uparrow}$. The proofs of the
following results can be found in Chapter III of Bertoin \cite{be}.
\begin{proposition}
$\mathrm{(i)}$ If $\mathrm{d}>0$ one has
\[
\lim_{x\to 0}\frac{U^{\uparrow}_x}{x}=\mathrm{d}\qquad \textrm{almost surely.}
\]
$\mathrm{(ii)}$ If $\mathrm{d}=0$ and $f:[0, \infty)\to [0, \infty)$ is an increasing function such that
$t\to f(t)/t$ increases, one has
\[
\limsup_{x\to 0}\frac{U^{\uparrow}_x}{f(x)}=\infty \quad
\textrm{a.s.} \quad\textrm{if and only if}\quad
\int_{0}\bar{\nu}(f(t))\ud t=\infty,
\]
where $\bar{\nu}(t)=\nu((t,\infty))$.\\
Moreover, 
\[
\textrm{if }\quad \int_{0}\bar{\nu}(f(t))\ud t<\infty\qquad\textrm{then}\qquad
\lim_{x\to 0}\frac{T^{\uparrow}_x}{f(x)}=\lim_{x\to 0}\frac{U^{\uparrow}_x}{f(x)}=0 \qquad \textrm{almost surely.}
\]
$\mathrm{(iii)}$ If $\e(T_1)<\infty$ one has
\[
\lim_{x\to =+\infty}\frac{U^{\uparrow}_x}{x}=\e(T_1)\qquad
\textrm{almost surely.}
\]
$\mathrm{(iv)}$ If $\e(T_1)$ is infinite and $f:[0, \infty)\to [0, \infty)$ is an increasing function such that the mapping
$t\to f(t)/t$ increases, one has
\[
\limsup_{x\to +\infty }\frac{U^{\uparrow}_x}{f(x)}=\infty\quad
\textrm{a.s.} \quad\textrm{if and only if}\quad
\int^{+\infty}\bar{\nu}(f(t))\ud t=\infty.
\]
Moreover, 
\[
\textrm{if }\quad \int^{+\infty}\bar{\nu}(f(t))\ud t<\infty\qquad\textrm{then}\qquad
\lim_{x\to +\infty}\frac{T^{\uparrow}_x}{f(x)}=\lim_{x\to +\infty}\frac{U^{\uparrow}_x}{f(x)}=0 \qquad\textrm{almost surely.}
\]
\end{proposition}
\section{Proofs of the main results.}
For simplicity, we introduce the notation
\[
h(t)=\frac{\log|\log t|}{\Phi(t^{-1}\log|\log t|)}.
\]
We start with the proof of the first part of Theorem 3, since a key
result on subordinators due to Fristed and Pruitt \cite{FP} easily
yields the result. The second part will be proved after the proof of Theorem 2,
since the latter is necessary for its proof.\\
\textit{Proof of Theorem 3 (first part):} First we will observe that
$\psi(\lambda)=O(\lambda^2)$, as $\lambda$ goes to $+\infty$, then
$\lambda^{1/2}=O(\Phi(\lambda))$. Since the last passage time
process $U^{\uparrow}$ is a subordinator with Laplace exponent
$\Phi$ and the future infimum process is the right-inverse of the
last passage times $U^{\uparrow}$, then according to Theorem 2 and
Remark on p. 176 in Fristed and Pruitt \cite{FP}, there exists a
positive constant $c$ such that
\[
\limsup_{t\to 0}\frac{J_t\Phi(t^{-1}\log|\log t|)}{\log|\log
t|}=c, \qquad \p^{\uparrow}-a.s.,
\]
and
\[
\limsup_{t\to +\infty}\frac{J_t\Phi(t^{-1}\log \log t)}{\log \log
t}=c, \qquad \p^{\uparrow}-a.s.,
\]
then the first part of Theorem 3 is proved.\QED
\noindent\textit{Proof of Theorem 2:} We only prove the result for
small times since the proof for large times is very similar. The
lower bound is easy to deduce from Theorem 3 and since
$J_t^{\uparrow}\leq \xi^{\uparrow}_t$, where $J_t^{\uparrow}$
denotes the future infimum of $\xi^{\uparrow}$. Hence
\[
c= \limsup_{t\to 0}\frac{J_t\Phi(t^{-1}\log|\log t|)}{\log|\log
t|}\leq \limsup_{t\to 0}\frac{\xi_t\Phi(t^{-1}\log|\log
t|)}{\log|\log t|}\qquad \p^{\uparrow}-a.s.
\]
Now, we prove the upper bound. Let $(x_n)$ be a decreasing sequence
of positive real numbers which converges to $0$, in particular we
choose $x_n=r^n$, for $r<1$.\\
Recall that $S$ is the supremum process of $\xi$, i.e. $S_t=sup_{0\leq u\leq t}\xi_u$. We define the events $A_n=\{S_{x_{n}}>\eta c h(x_{n+1})\}$, where
$\eta \geq c^{-1}(2+r^{-1})$ and $S$ is the supremum process. From the first Borel-Cantelli's Lemma,
if $\sum_{n}\p^{\uparrow}(A_n)<\infty$, it follows
\[
S_{r^{n}}\leq \eta c h(r^{n+1})\qquad\p^{\uparrow}\textrm{-a.s.,}
\]
for all large $n$. Since the function $h$ and the process $S$ are
increasing in a neighbourhood of 0, we have
\[
S_t\leq \eta c h(t)\qquad\textrm{for }\quad r^{n+1}\leq t\leq
r^n,\quad \textrm{under }\p^{\uparrow}.
\]
Then, it is enough to prove that
$\sum_{n}\p^{\uparrow}(A_n)<\infty$. In this direction, we will
prove the following lemma,
\begin{lemma} Let $0<\epsilon<1$ and $r<1$. If we assume that condition (H2) is satisfied
then there exists a positive constant $C(\epsilon)$ such that
\begin{equation}\label{in2}
\p^{\uparrow}\big(J_{r^n}>(1-\epsilon)\eta ch(r^{n+1})\big) \geq
C(\epsilon) \p^{\uparrow}\big(A_n\big) \quad\textrm{as }\quad n\to
+\infty.
\end{equation}
\end{lemma}
\textit{Proof of Lemma 3:} From the inequality (\ref{deslogr}), we have that
\[
\p^{\uparrow}\big(J_{r^n}>(1-\epsilon)\eta
ch(r^{n+1})\big)\geq\p^{\uparrow}_{\eta
ch(r^{n+1})}\big(J_{0}>(1-\epsilon)\eta ch(r^{n+1})\big)
\p^{\uparrow}\big(S_{r^n}>\eta ch(r^{n+1}) \big),
\]
and since $J_{0}$ is the absolute minimum of $(\xi,
\p^{\uparrow}_{\eta ch(r^{n+1})})$ then by (\ref{distabsmin})
\[
\p^{\uparrow}_{\eta ch(r^{n+1})}\big(J_{0}>(1-\epsilon)\eta
ch(r^{n+1})\big) =\frac{W\big(\epsilon \eta ch(r^{n+1})\big)}{W\big(
\eta ch(r^{n+1})\big)}.
\]
On the other hand, an application  of Proposition III.1 in Bertoin
\cite{be} gives that there exist a positive real number $K_1$ such
that
\begin{equation}\label{relwpsi}
 K_1\frac{1}{x\psi(1/x)}\leq W(x)\leq K_1^{-1}\frac{1}{x\psi(1/x)},
 \qquad\textrm{ for all} \quad x>0,
\end{equation}
then it is clear that
\[
\frac{W\big(\epsilon \eta ch(r^{n+1})\big)}{W\big( \eta
ch(r^{n+1})\big)}\geq K_1^2\epsilon^{-1}\frac{\psi\big(1/\eta
ch(r^{n+1})\big)}{\psi\big(\epsilon^{-1}/\eta ch(r^{n+1})\big)}.
\]
From this inequality and  condition (H2), there exist a
positive constant $C(\epsilon)$ such that  for $n$ sufficiently
large
\[
\p^{\uparrow}\big(J_{r^n}>(1-\epsilon)\eta ch(r^{n+1})\big) \geq
C(\epsilon) \p^{\uparrow}\big(S_{r^n}>\eta ch(r^{n+1})\big),
\]
which proves our result.\QED \noindent Now, we  prove the upper
bound for the law of the iterated logarithm of $(\xi,
\p^{\uparrow})$. Fix $0<\epsilon<1/(2+r^{-1})$. Since $J$ can be
seen as the right inverse of $U$, it is straightforward that
\[
\p^{\uparrow}\big(J_{r^{n}}>(1-\epsilon)\eta ch(r^{n+1})
\big)=\p^{\uparrow}\big(U_{(1-\epsilon)\eta ch(r^{n+1})}<r^n \big),
\]
and this probability is bounded from above by
\[
\exp\{\lambda r^n\}\e^{\uparrow}\Big(\exp\big\{-\lambda
U_{(1-\epsilon)\eta c h(r^{n+1})} \big\}\Big)=\exp\big\{\lambda r^n
-(1-\epsilon)\eta c h(r^{n+1})\Phi(\lambda)\big\},
\]
for every $\lambda\geq 0$. We choose $\lambda=r^{-(n+1)}\log|\log
r^{n+1}|$, then
\[
\p^{\uparrow}\big(J_{r^{n}}>(1-\epsilon)\eta ch(r^{n+1})
\big)\leq\exp\Big\{-\big((1-\epsilon)\eta c-r^{-1}\big)\log|\log
r^{n+1}|\Big\},
\]
hence from the above inequality and Lemma 3, we have that
\[
C(\epsilon)\sum_{n}\p^{\uparrow}(A_n)\leq C_1\sum_{n}\big(\log
(n+1) \big)^{(1-\epsilon)\eta c -r^{-1}}<+\infty,
\]
since $(1-\epsilon)\eta c -r^{-1}>1$.\\
Hence, we have
\[
\limsup_{t\to 0}\frac{S_t}{h(t)}\leq \eta c,
\qquad\p^{\uparrow}\textrm{-a.s.,}
\]
for $\eta c>3$, since we can choose $r$ close enough to $1$.\\
The two preceding parts show that
\[
\limsup_{t\to 0}\frac{\xi_t}{h(t)}\in [c,\eta c],
\qquad\p^{\uparrow}\textrm{-a.s.}
\]
By the Blumenthal zero-one law, it must be a constant number $k$,
$\p^{\uparrow}-$a.s.\QED \noindent\textit{Proof of Theorem 3 (second
part):} First we prove the result for large times. Assume that the
additional hypothesis (H2) is satisfied. Since $J_t\leq
\xi_t$ for every $t\geq 0$ and Theorem 2, it is clear that
\[
\limsup_{t\to +\infty}\frac{J_t}{h(t)}\leq \limsup_{t\to
+\infty}\frac{\xi_t}{h(t)}=k \qquad\p^{\uparrow}\textrm{-a.s.},
\]
then the upper bound is proved.\\
Now, fix $\epsilon \in (0, 1/2)$ and define
\[
R_n=\inf\left\{s\geq n:\frac{\xi^{\uparrow}_s}{kh(s)}\geq
(1-\epsilon)\right\}.
\]
From the above definition, it is clear that $R_n\geq n$ and that
$R_n$ diverge a.s. as $n$ goes to $+\infty$. From Theorem 2, we
deduce that $R_n$ is finite a.s.\\
Now, by (\ref{distabsmin}) and since $(\xi, \p^{\uparrow})$ is a
strong Markov process with no positive jumps, we have that
\[
\begin{split}
\p^{\uparrow}\left(\frac{J_{R_n}}{k h(R_n)}\geq
(1-2\epsilon)\right)&=\p^{\uparrow}\left(J_{R_n}\geq
\frac{(1-2\epsilon)\xi_{R_n}}{(1-\epsilon)}\right)\\
&=\e^{\uparrow}\bigg(\p^{\uparrow}\left(J_{R_n}\geq
\frac{(1-2\epsilon)\xi_{R_n}}{(1-\epsilon)}\Big|\xi_{R_n}\right)\bigg)\\
&=\e^{\uparrow}\bigg(\frac{W(\epsilon
\xi_{R_n})}{W(\xi_{R_n})}\bigg).
\end{split}
\]
Now applying (\ref{relwpsi}), we have that
\[
\e^{\uparrow}\bigg(\frac{W(\epsilon
\xi_{R_n})}{W(\xi_{R_n})}\bigg)\geq K_1^2
\epsilon^{-1}\e^{\uparrow}\bigg(\frac{\psi(1/\xi_{R_n})}{\psi(\epsilon^{-1}/\xi_{R_n})}\bigg)
\]
and since the  hypothesis (H2) is satisfied,
an application of the Fatou-Lebesgue Theorem shows that
\[
\liminf_{n\to
+\infty}\e^{\uparrow}\bigg(\frac{\psi(1/\xi_{R_n})}{\psi(\epsilon^{-1}/\xi_{R_n})}\bigg)\geq
\e^{\uparrow}\bigg(\liminf_{n\to
+\infty}\frac{\psi(1/\xi_{R_n})}{\psi(\epsilon^{-1}/\xi_{R_n})}\bigg)
>0,
\]
which implies that
\[
\lim_{n\to +\infty}\p^{\uparrow}\left(\frac{J_{R_n}}{k h(R_n)}\geq
(1-2\epsilon)\right)>0.
\]
Since $R_n\geq n$,
\[
\p^{\uparrow}\left(\frac{J_t}{k h(t)}\geq(1-2\epsilon), \textrm{ for
some } t\geq n\right)\geq\p^{\uparrow}\left(\frac{J_{R_n}}{k
h(R_n)}\geq (1-2\epsilon)\right).
\]
Therefore, for all $\epsilon \in (0, 1/2)$
\[
\p^{\uparrow}\left(\frac{J_t}{k h(t)}\geq(1-2\epsilon), \textrm{
i.o., as } t\to +\infty\right)\geq\lim_{n\to
+\infty}\p^{\uparrow}\left(\frac{J_{R_n}}{k h(R_n)}\geq
(1-2\epsilon)\right)>0.
\]
The event on the left hand side is in the upper-tail sigma-field
$\cap_{t}\sigma\{\xi^{\uparrow}_s: s\geq t\}$ which is trivial from
Bertoin's construction of $(\xi, \p^{\uparrow})$ (see Theorem VII.20
in \cite{be}). Hence
\[
\limsup_{t\to +\infty}\frac{J_t}{ h(t)}\geq k(1-2\epsilon), \qquad
\p^{\uparrow}-\textrm{a.s.},
\]
and since $\epsilon$ can be chosen arbitrarily small, the result
for large times is proved.\\
In order to prove the  law of the iterated logarithm for small
times, we now define the following stopping time
\[
R_n=\inf\left\{\frac{1}{n}<s:\frac{\xi^{\uparrow}_s}{kh(s)}\geq (1-\epsilon)\right\}.
\]
Following same argument as above and assuming that (H1) is satisfied, we get that for a fixed $\epsilon \in (0,1/2)$ and $n$ sufficiently large
\[
\p^{\uparrow}\left(\frac{J_{R_n}}{k h(R_n)}\geq
(1-2\epsilon)\right)>0.
\]
Next, we note that
\[
\p^{\uparrow}\left(\frac{J_{R_p}}{k h(R_p)}\geq
(1-2\epsilon), \textrm{ for some }p\geq n\right)\geq \p^{\uparrow}\left(\frac{J_{R_n}}{k h(R_n)}\geq
(1-2\epsilon)\right).
\]
Since $R_n$ converge a.s. to $0$ as $n$ goes to $\infty$, the conclusion follows taking the limit when $n$ goes towards to $+\infty$.\QED
\noindent\textit{Proof of Theorem 4:} The proof of this theorem is
very similar to the proof of the previous result. Following  same
arguments, we first prove the law of the iterated logarithm for
large times. Assume that the hypothesis (H2) and (H3)
are satisfied. Since $\xi^{\uparrow}_t -J^{\uparrow}_t\leq
\xi^{\uparrow}_t$ for every $t$ and Theorem 2, it is clear that
\[
\limsup_{t\to +\infty}\frac{\xi_t -J_t}{h(t)}\leq \limsup_{t\to
+\infty}\frac{\xi_t}{h(t)}=k \qquad\p^{\uparrow}\textrm{-a.s.},
\]
then the upper bound is proved.\\
Now, fix $\epsilon \in (0, 1/2)$ and similarly as the last proof
we define
\[
R_n=\inf\left\{s\geq n:\frac{\xi^{\uparrow}_s}{kh(s)}\geq
(1-\epsilon)\right\}.
\]
Now, by (\ref{distabsmin}) and since $(\xi, \p^{\uparrow})$ is a
strong Markov process with no positive jumps, we have that
\[
\begin{split}
\p^{\uparrow}\left(\frac{\xi_{R_n}-J_{R_n}}{k h(R_n)}\geq
(1-2\epsilon)\right)&=\p^{\uparrow}\left(J_{R_n}\leq
\frac{\epsilon}{(1-\epsilon)}\xi_{R_n}\right)\\
&=\e^{\uparrow}\bigg(\p^{\uparrow}\left(J_{R_n}\leq
\frac{\epsilon}{(1-\epsilon)}\xi_{R_n}\Big|\xi_{R_n}\right)\bigg)\\
&=1-\e^{\uparrow}\bigg(\frac{W\big(k(\epsilon)\xi_{R_n}
\big)}{W(\xi_{R_n})}\bigg),
\end{split}
\]
where $k(\epsilon)=(1-2\epsilon)/(1-\epsilon)$.\\
Since the  hypothesis (H3) is satisfied, an
application of the Fatou-Lebesgue  Theorem shows that
\[
\limsup_{n\to +\infty}\e^{\uparrow}\bigg(\frac{W\big(k(\epsilon)
\xi_{R_n}\big)}{W(\xi_{R_n})}\bigg)\leq\e^{\uparrow}\bigg(\limsup_{n\to
+\infty}\frac{W\big(k(\epsilon)
\xi_{R_n}\big)}{W(\xi_{R_n})}\bigg)<1,
\]
which implies that
\[
\lim_{n\to +\infty}\p^{\uparrow}\left(\frac{\xi_{R_n} -J_{R_n}}{k
h(R_n)}\geq (1-2\epsilon)\right)>0.
\]
Again, since $R_n\geq n$,
\[
\p^{\uparrow}\left(\frac{\xi_t -J_t}{k h(t)}\geq(1-2\epsilon),
\textrm{ for some } t\geq
n\right)\geq\p^{\uparrow}\left(\frac{\xi_{R_n} -J_{R_n}}{k
h(R_n)}\geq (1-2\epsilon)\right).
\]
Therefore, for all $\epsilon \in (0, 1/2)$
\[
\p^{\uparrow}\left(\frac{\xi_t -J_t}{k h(t)}\geq(1-2\epsilon),
\textrm{ i.o., as } t\to +\infty\right)\geq\lim_{n\to
+\infty}\p^{\uparrow}\left(\frac{\xi_{R_n} -J_{R_n}}{k h(R_n)}\geq
(1-2\epsilon)\right)>0.
\]
The event on the left hand side is in the upper-tail sigma-field
$\cap_{t}\sigma\{\xi^{\uparrow}_s: s\geq t\}$ which is trivial, then
\[
\limsup_{t\to +\infty}\frac{\xi_t -J_t}{ h(t)}\geq k(1-2\epsilon),
\qquad \p^{\uparrow}-\textrm{a.s.},
\]
and since $\epsilon$ can be chosen arbitrarily small, the result
for large times is proved.\\
Similarly as in the proof of the previous result, we can prove the
result for small times using the following stopping time
\[
R_n=\inf\left\{\frac{1}{n}<s:\frac{\xi^{\uparrow}_s}{k
h(s)}\geq (1-\epsilon)\right\}.
\]
Following same argument as above and assuming that (H1) and (H4) are satisfied, we get that for a fixed $\epsilon \in (0,1/2)$ and $n$ sufficiently large
\[
\p^{\uparrow}\left(\frac{\xi_{R_n}-J_{R_n}}{k h(R_n)}\geq
(1-2\epsilon)\right)>0.
\]
Next, we note that
\[
\p^{\uparrow}\left(\frac{\xi_{R_n}-J_{R_p}}{k h(R_p)}\geq
(1-2\epsilon), \textrm{ for some }p\geq n\right)\geq \p^{\uparrow}\left(\frac{\xi_{R_n}-J_{R_n}}{k h(R_n)}\geq
(1-2\epsilon)\right).
\]
Again, since $R_n$ converge a.s. to $0$ as $n$ goes to $\infty$, the conclusion follows taking the limit when $n$ goes towards to $+\infty$.\QED
\noindent \textit{Proof of Theorem 5:} Let $(x_{n})$ be a decreasing
sequence such that $\lim x_{n}=0$. We define the events
\[
A_{n}=\Big\{\textrm{ There exist }t\in [U^{\uparrow}_{x_{n+1}},
U^{\uparrow}_{x_{n}}]\textrm{ such that }
\xi^{\uparrow}_{t}<f(t)\Big\}.
\]
Since $U^{\uparrow}_{x_{n}}$ tends to 0, a.s. when $n$ goes to
$+\infty$ , we have
\[
\Big\{ \xi^{\uparrow}_{t}<f(t), \textrm{ i.o., as } t\to
0\Big\}=\limsup_{n\to +\infty}A_{n}.
\]
Let us chose $x_{n}=r^{n}$, for $r<1$. Since $f$ is increasing the
following inclusions hold
\[
A_{n}\subset \Big\{ \textrm{There exist } t\in [r^{n+1}, r^{n}]
\textrm{ such that }tr<f\big(U^{\uparrow}_{t}\big)\Big\},
\]
and
\[
\Big\{ \textrm{There exist } t\in [r^{n+1}, r^{n}] \textrm{ such
that }tr^{-1}<f\big(U^{\uparrow}_{t}\big)\Big\}\subset A_{n}
\]
Then we prove the convergent part. Let us suppose that $f$ satisfies
\[
\int_{0}f(x)\nu(\ud x)<\infty.
\]
Hence from Theorem VI.3.2 in \cite{GS} and the fact that
$U^{\uparrow}$ is a subordinator, we have that
\[
\p^{\uparrow}\Big( tr<f\big(U_{t}\big), \textrm{ i.o., as } t\to
0\Big)=0,
\]
which implies that
\[
\lim_{t\to 0}\frac{\xi_{t}}{f(t)}=\infty \qquad
\p^{\uparrow}\textrm{-a.s.},
\]
since we can replace $f$ by $cf$, for any $c>1$.\\
Similarly, if $f$ satisfies that
\[
\int_{0}f(x)\nu(\ud x)=\infty,
\]
again from Theorem VI.3.2 in \cite{GS}, we have that
\[
\p^{\uparrow}\Big( tr^{-1}<f\big(U_{t}\big), \textrm{ i.o., as }
t\to 0\Big)=1,
\]
which implies that
\[
\liminf_{t\to 0}\frac{\xi_{t}}{f(t)}=0\qquad
\p^{\uparrow}\textrm{-a.s.},
\]
since we can replace $f$ by $cf$, for any $c<1$.\\
The integral test at $+\infty$ is very similar to this of small
times, it is enough to take $x_{n}=r^{n}$, for $r>1$ and follows the
same arguments as in the proof for small times.\QED
\noindent\textit{Proof of Theorem 6:} The proof of parts $(ii)$ and
$(iv)$ follows from the proof of Theorem 5, it is enough to note
that we can replace $\xi^{\uparrow}$ by $J^{\uparrow}$ in the sets
$A_{n}$. The proof of parts $(i)$ and $(iii)$ follows from
Proposition 4.4 in \cite{be5}.\QED \noindent\textit{Proof of
Proposition 1:} Similarly as in the proof of Theorem 5, let
$(x_{n})$ be a decreasing sequence such that $\lim x_{n}=0$ and
$c>1$. We define the events
\[
A_{n}=\Big\{\textrm{ There exist }t\in [U^{\uparrow}_{x_{n+1}},
U^{\uparrow}_{x_{n}}]\textrm{ such that }
\xi^{\uparrow}_{t}<cf(t)\Big\}.
\]
Since $U^{\uparrow}_{x_{n}}$ tends to 0, a.s. when $n$ goes to
$+\infty$ , we have
\[
\Big\{ \xi^{\uparrow}_{t}<cf(t), \textrm{ i.o., as } t\to
0\Big\}=\limsup_{n\to +\infty}A_{n}.
\]
Since $f$ is increasing the following inclusion holds
\[
A_{n}\subset \Big\{x_{n+1}<cf\big(U^{\uparrow}_{x_{n}}\big) \Big\}.
\]
On the other hand
\[
\p^{\uparrow}\Big(x_{n+1}<cf\big(U_{x_{n}}\big)
\Big)=\p^{\uparrow}\Big(f^{-1}\big(x_{n+1}/c\big)<U_{x_{n}} \Big),
\]
where $f^{-1}$ is the right-inverse of $f$.\\
Now, we take $x_{n}=cf(r^n)$, for $r<1$. Since $f$ is increasing and
from the above equality,  we get that
\[
\p^{\uparrow}\Big(x_{n+1}<cf\big(U_{x_{n}}\big) \Big)\leq
\p^{\uparrow}\Big(r^{n+1}<U_{cf(r^{n})} \Big)
\]
The obvious inequality
\[
\p^{\uparrow}\Big(a<U_{t} \Big)\leq
(1-e^{-1})^{-1}\Big(1-\exp\big\{-t\Phi(1/a)\big\}\Big),
\]
applied for $t=cf(r^{n})$ and $a=r^{n+1}$ entails that
\[
\p^{\uparrow}\Big(r^{n+1}<U_{cf(r^{n})} \Big)\leq (1-e^{-1})^{-1}c
f(r^{n})\Phi(r^{-(n+1)}).
\]
Since the mapping $t\to t\Phi(1/t)$ increases, it is not difficult to deduce that the function $\Phi$ satisfies that
\[
\Phi(r^{-(n+1)})\leq r^{-2}\Phi(r^{(n-1)}).
\]
Hence,
\[
\begin{split}
\sum_{n\geq k} \p^{\uparrow}\Big(r^{n+1}<U_{cf(r^{n})} \Big)&\leq
C(r)\sum_{n\geq k}\int_{n-1}^{n} f(r^{t})\Phi(r^{-t})\ud
t\\
&\leq C(r) \int_{0}^{r^{k-1}} x^{-1}f(x)\Phi(1/x)\ud x.
\end{split}
\]
Since the last integral is finite, by the Borel-Cantelli lemma, we
deduce that
\[
\p^{\uparrow} \Big( \xi_{t}<cf(t), \textrm{ i.o., as } t\to
0\Big)=0,
\]
for all $c\geq 1$, hence
\[
\lim_{t\to 0}\frac{\xi_{t}}{f(t)}=\infty, \qquad
\p^{\uparrow}\textrm{-a.s.}
\]
In order to prove that the future infimum satisfies the same result,
 we note first that we can replace $\xi^{\uparrow}$ by its future infimum in
the sets $A_{n}$, and then the same arguments will give us the
desired result.\QED

\noindent{\bf Acknowledgements.} I would like to thank Lo\"\i c
Chaumont for guiding me through the development of this work, and
for all his helpful advice. I also like to express my gratitude to Andreas Kyprianou for all his useful suggestions.

\end{document}